\documentclass[english,10pt]{amsart}
\usepackage[francais]{babel}
\usepackage{amsmath,amsfonts,amssymb,multicol, amsthm}
\usepackage{amsrefs}
\usepackage{mathrsfs}
\usepackage[utf8]{inputenc}
\usepackage{esint}
\usepackage{color, framed}

\usepackage{graphicx}
\usepackage[normalem]{ulem}

\newcommand{\R}{\ensuremath{\mathbb{R}}}

\newcommand{\N}{\ensuremath{\mathbb{N}}}

\renewcommand{\leq}{\leqslant}
\renewcommand{\geq}{\geqslant}

\renewcommand\proof{\noindent {\bf Proof.} \quad}
\newcommand\eproof{{\hfill $\square$}}

\newtheorem{theorem}{Theorem}
\newtheorem{thmx}{Theorem}

\newtheorem{lemma}{Lemma}

\theoremstyle{remark}
\newtheorem{remark}[theorem]{Remark}

\newcommand{\super}{\overline}

\newcounter{numeroexo}

\title{Classification and Liouville-type theorems for semilinear elliptic equations in unbounded domains.}
\author{Louis Dupaigne}
\address{Institut Camille Jordan, UMR CNRS 5208, Universit\'e Claude Bernard Lyon 1, 43 boulevard du 11 novembre 1918, 69622 Villeurbanne cedex, France}
\email{dupaigne@math.univ-lyon1.fr}
\author{Alberto Farina}
\address{LAMFA, UMR CNRS 7352, Universit\'e Picardie Jules Verne 33, rue St Leu, 80039 Amiens, France}
\email{alberto.farina@u-picardie.fr}
\date{}

\begin{document}

\maketitle

\section{Introduction and main results} 

A noncompact Riemannian manifold $\mathcal M$ is said to be parabolic if every positive superharmonic function $u:\mathcal M\to\R_+$ is constant. This is equivalent to asking that there exists no positive fundamental solution of the Laplace equation or that Brownian motion is recurrent on $\mathcal M$, see e.g. \cite{gri}. $\R^2$ is parabolic, while $\R^N$ is not when $N\ge 3$. Hoping that our findings could give an interesting way to study noncompact Riemannian manifolds\footnote{at least under curvature or curvature-dimension conditions.} in higher dimensions, we focus in what follows on those positive superharmonic functions that are also solutions of a semilinear elliptic equation. We shall be interested more specifically in solutions having some stability property, although this assumption is not needed for some of our results. We work in the Euclidean setting and prove sharp Liouville-type theorems and classification results.

\medskip

Since $\R^N$ is invariant under dilations, such theorems are very much related to the corresponding regularity theory, notably some tools and results of the recent paper \cite{cfrs} exploiting the dilation invariance (a version of Pohozaev's classical argument \cite{poho} suited for stable solutions), the autonomous nature of the equation (a geometric Poincaré formula\footnote{which also holds on manifolds, see \cite{FMV}, and ca also be reformulated in the general language of $\Gamma$-calculus, see \cite{v}} discovered in \cite{sz}) and universal bounds in the $C^\alpha$ and $H^1$ norms (due to \cite{cfrs}). Among novelties in what follows, we are able to handle the critical dimension $N=10$, to work without any convexity assumption on the nonlinearity (except for the regularity theory of finite Morse index solutions) and to classify solutions which are stable outside a compact set. Finally, we provide new Liouville-type theorems on half-spaces, improving those of \cite{FVarma}, as well as on coercive epigraphs.

\medskip

Our first result extends a theorem in our earlier work \cite{DuFa}, which held for bounded stable solutions in dimension $ N \leq 4$, and reads as follows.
\begin{theorem} \label{theorem:liouville0} Assume that $u\in C^2(\R^N)$ is bounded below and that $u$ is a stable solution of 
\begin{equation}\label{equazione}
-\Delta u = f(u)\qquad\text{in $\R^N$}
\end{equation}
where $f:\R\to\R$ is locally Lipschitz and nonnegative. If $N\le 10$, then $u$ must be constant.
\end{theorem}
We recall that if $f\in C^1(\R)$, a solution $u$ of $-\Delta u = f(u)$ in $ \Omega \subseteq \R^N$  is stable if for every $\varphi\in C^1_c(\Omega)$, there holds
\begin{equation}\label{stab}
\int_{\Omega}f'(u)\varphi^2 \le \int_{\Omega}\vert\nabla\varphi\vert^2
\end{equation}
When $f$ is locally Lipschitz, the definition\footnote{In \cite{cfrs}, a slightly different definition than that in \cite{fsv}, is used. Note that when $f$ is nondecreasing, the condition in \cite{cfrs} is {\it a priori} more restrictive. However, as observed by the authors of both papers, when restricting to test functions of the form $\varphi =\vert\nabla u\vert\psi$ or of the form $\varphi=\partial_r u\psi$, $\psi\in C^1_c(\R^N)$, then both definitions yield the same information. In particular, our results hold using one definition or the other.} requires more care, see \cite{fsv} and \cite{cfrs}. 

\medskip

Under a weaker lower bound, the previous result remains true up to dimension 9:

\begin{theorem} \label{theorem:liouville0b} Assume that $u\in C^2(\R^N)$  is a stable solution of \eqref{equazione},
where $f:\R\to\R$ is locally Lipschitz and nonnegative. Assume in addition that  for some $C>0$
\begin{equation}\label{log}u(x)\ge -C\ln(2+\vert x\vert), \qquad x\in\R^N.
\end{equation}
If $N\le 9$, then $u$ must be constant.
\end{theorem}

\begin{remark}\label{Rem1}

\

Some remarks are in order.

\begin{itemize}
\item Theorem \ref{theorem:liouville0} is sharp. Indeed,  if $N\ge11$, for $f(u)=u^p$, $p$ sufficiently large, there exists a nontrivial positive bounded stable solution to the equation, see \cite{fa2}.
\item Theorem \ref{theorem:liouville0b} is also sharp. Indeed, if $N\ge 10$ and $f(u)=e^u$, there exists\footnote{To prove this, using Emden's transformation, the equation is equivalent to an autonomous ode having a unique stationary point, which corresponds to the singular solution $u_s(x)=\ln\frac{2(N-2)}{\vert x\vert^2}$ and is attractive, see e.g. pp. 36-37 in Chapter 2 in \cite{dup}. In addition, $u_s$ is stable for $N\ge10$, thanks to Hardy's inequality. Since $u_s$ is singular at the origin, it is clear that $u\le u_s$ in some ball $B_R$. In fact, the inequality holds throughout $\R^N$ and so $u$ is also stable. Otherwise, there would exist $R'$ such that $u_s-u=0$ on $\partial B_{R'}$. Using $\varphi=u_s-u\in H^1_0(B_{R'})$ as a test function in \eqref{stab} would then contradict the stability of $u_s$, using (the proof of) Proposition 3.2.1 in \cite{dup}. } a radial stable solution $u\in C^2(\R^N)$ such that $u(x)\sim -2\ln\vert x\vert$ as $\vert x\vert\to+\infty$.

\item In the important cases where $f(u) = \vert u \vert^{p-1}u$, $p>1$ and $f(u) = e^{u} $, Theorems \ref{theorem:liouville0} and \ref{theorem:liouville0b} were already known to hold, see \cites{fa2,FaEXP}. In addition, in these cases, the theorems hold without assuming any bound on $u$, due to the scale invariance of the equations. 

\item In the particular case where $u$ is radial and bounded, Theorems \ref{theorem:liouville0} and \ref{theorem:liouville0b} were already known to hold, see \cites{cc, villegas}.
\item It will be clear from the proof of Theorem \ref{theorem:liouville0b} that the same conclusion holds true if we replace \eqref{log} by the weaker one : $ u(x)\ge -C\ln^{\gamma}(2+\vert x\vert)$, where $\gamma\ge1$.
Some lower bound is however needed {in Theorems \ref{theorem:liouville0} and \ref{theorem:liouville0b}.} Indeed, the function $u(x)=-\vert x\vert^2$ is a stable solution (with $f(u)=2N$). It is bounded above, but not below. 
\item The assumption that $u$ is superharmonic (i.e. $f\ge0$) is also essential. Indeed, if $f(u)=u-u^3$, then $f$ changes sign and for any $N\ge1$, $u(x)=\tanh(x_1/\sqrt 2)$ is bounded  and monotone (hence stable) yet non constant\footnote{One could ask whether other kinds of solutions exist. This is a delicate question deeply related to the celebrated De Giorgi conjecture, see e.g. \cite{fvdg,cw,dup} for more on this subject.}. Similarly, if $f(u)=-2N$, then $f\le 0$ and $u(x)=\vert x\vert^2$ is a stable solution bounded below. Still, when $f\le0$, one can clearly use our theorems to classify stable solutions which are bounded {\it above} (simply work with $-u$ in place of $u$).
\item At last, the theorems cannot be generalized to solutions which are merely stable outside a compact set (and so to finite Morse index solutions). Indeed, if $N\ge 3$, the standard bubble is a nonconstant bounded solution of $-\Delta u = u^{\frac{N+2}{N-2}}$ in $\R^N$ of Morse index 1. The same remark holds in dimension $N=2$ for solutions stable outside a compact set satisfying the weaker bound \eqref{log}, by considering the Liouville equation $-\Delta u =e^u$ in $\R^2$ (see Theorem 3 in \cite{FaEXP}).
\item Additional examples of non-trivial solutions of finite Morse index are provided by {\it subcritical} nonlinearities of the form $f(u)=[(u-\beta)^+]^p$, for all $p\in(1,\frac{N+2}{N-2})$, $N\ge 3$ and some $\beta>0$, see \cite{df2}.
\item Note however that all these counter-examples are radial functions. We address the question of radial symmetry in Theorem \ref{theorem:symmetry}  below.
\end{itemize}
\end{remark}

As important corollaries, we obtain the following Liouville-type results on half-spaces and on coercive epigraphs.
\begin{theorem}\label{th:liouville1}
Let $u\in C^2(\R^N_+)$ be a bounded solution of 
\begin{equation}\label{pbl-half-space}
\left\{
\begin{aligned}
-\Delta u& = f(u)\quad\text{in $\R^N_+$,}\\
u&>0\quad\text{in $\R^N_+$,}\\
u&=0\quad\text{on $\partial \R^N_+.$}
\end{aligned}
\right.
\end{equation}
Assume $f\in C^1(\R)$ and 
\begin{enumerate}
	\item either $ f(t) \geq 0 $ for $ t \geq 0,$
	\item or there exists $z>0$ such that $ f(t) \geq 0$ for $ t \in [0,z]$ and $ f(t) \leq 0$ for $t \geq z$.
\end{enumerate}
If $ 2 \leq N\le 11$,  then $u$ must be one-dimensional and monotone (i.e., $u=u(x_N)$ and $ \frac{\partial u }{\partial x_N} >0$ in 
$ \R^N_+$). 
\end{theorem}

\medskip

Theorem \ref{th:liouville1} recovers and improves upon a result of \cite{FVarma}, which held for $N \leq 5$ (see also \cite{BCN, BCNMagenes, CLZ, Dancer, dancer2, Fa, fa2, Fa3, FaSciunzi, FaSciunzi2, fsv} for some other results concerning problem \eqref{pbl-half-space}.)

\medskip

For the Neumann boundary condition, the following result holds.
\begin{theorem}\label{th:liouville2}
Let $u\in C^2(\R^N_+)$ be a stable solution of 
\begin{equation}\label{Neumann}
\left\{
\begin{aligned}
-\Delta u& = f(u)\quad\text{in $\R^N_+$,}\\
\partial_n u&=0\quad\text{on $\partial \R^N_+.$}
\end{aligned}
\right.
\end{equation}
Assume in addition that $f:\R\to\R$ is locally Lipschitz and nonnegative and $u$ is bounded below.
If $N\le 10$, then $u$ must be constant. 
\end{theorem}
Here, a solution of \eqref{Neumann} is said to be stable if \eqref{stab} holds with $\Omega = \R^N_+$ for all $\varphi\in C^1_c(\super{\R^N_+})$.

\medskip

At last, for coercive epigraphs, the following result holds true.
\begin{theorem}\label{th:liouville3}
Let $\Omega\subset\R^N$ denote a locally Lipschitz coercive epigraph and $u\in C^2(\Omega)\cap C(\super\Omega)$ be a bounded solution of 
\begin{equation}\label{pbl-epi}
\left\{
\begin{aligned}
-\Delta u& = f(u)\quad\text{in $\Omega$,}\\
u&\ge0\quad\text{in $\Omega$,}\\
u&=0\quad\text{on $\partial \Omega.$}
\end{aligned}
\right.
\end{equation}
Assume that $f\in C^1(\R)$, $f(t)>0$ for $t>0$ and $ 2 \leq N\le 11$. Then, $f(0)=0$ and $u=0$.
\end{theorem}

%
%

\noindent Before stating our results on solutions which are merely stable outside a compact set, it will be useful to discuss the proof of Theorem \ref{theorem:liouville0b}, which relies on the following {\it a priori} estimate, recently established in \cite{cfrs}. 

\begin{thmx}[\cite{cfrs}]\label{th:cfrs} Let $B_1$ be the unit ball of $\mathbb{R}^N$, $ N \geq 1$. Assume that $u\in C^2(B_1)$ is a stable solution of 
	$$
	-\Delta u = f(u)\quad\text{in $B_1$},
	$$
	where $f:\R\to\R$ is locally Lipschitz and nonnegative. If $N\le 9$, then 
	\begin{equation}\label{cfrs}
	\Vert u\Vert_{C^\alpha(\super{B_{1/2}})} \le C \Vert u\Vert_{L^1(B_{1})},
	\end{equation}
	where $\alpha\in(0,1)$, $C>0$ are dimensional constants.
\end{thmx}

Since the proof of Theorem \ref{theorem:liouville0b} is very short, we provide it without further ado.

\noindent {\bf Proof of Theorem \ref{theorem:liouville0b}. }

Given $R>2$, apply \eqref{cfrs} to $u_R(x)=u(Rx)$, leading to
$$
|u(x)-u(y)| \le C R^{-\alpha} |x-y|^{\alpha}\fint_{B_R}\vert u\vert \qquad \text{ for }x,y\in B_{R/2}.
$$
Since $u$ satisfies the lower bound \eqref{log}, observe that
$$
\vert u\vert = \vert u +C\ln(2+\vert x\vert) - C\ln(2+\vert x\vert)\vert\le \vert u+ C\ln(2+\vert x\vert)\vert + \vert C\ln(2+\vert x\vert)\vert = u+2C\ln(2+\vert x\vert).
$$
So, recalling that $u$ is superharmonic,
$$
0\le \fint_{B_R}\vert u\vert \le \fint_{B_R}u + 2C\ln (2+R)\le u(0)+2C\ln (2+R)
$$
and so 
$$
|u(x)-u(y)| \le C R^{-\alpha} |x-y|^{\alpha} (u(0)+2C\ln (2+R)).
$$
Let $R\to+\infty$ to conclude that $u(x)=u(y)$ for all $x,y\in\R^N$. 
\hfill\qed

\

\begin{remark}
Since Theorem \ref{theorem:liouville0b} fails for finite Morse index solutions, it follows from the proof above that the {\it a priori} estimate in Theorem \ref{th:cfrs} cannot hold either for such solutions.
\end{remark}

However, the equation is still smoothing, at least under an extra convexity assumption: 

\begin{theorem}
Let $B_1$ denote the unit ball in $\R^N$. Assume that $u\in H^1(B_1)$ is a weak solution of 
$$
-\Delta u = f(u)\qquad\text{in $B_1$}
$$
that has finite Morse index, where $f\ge0$ is nondecreasing and convex (hence locally Lipschitz) and $1\le N\le 9$. Then, $u\in C^{2,\alpha}(\super{B_{1/2}})$ for all $\alpha\in(0,1)$.
\end{theorem}

\proof
Indeed, according to Proposition 2.1. in \cite{ddf}, for every $x_0\in B_{1/2}$, there exists a ball $B(x_0,r_0)$ such that $u$ is stable in $B(x_0,r_0)$. In addition, according to Proposition 4.2 in \cite{cfrs}, there exists a sequence of $C^2$ stable solutions $(u_n)$ in $B(x_0,r_0)$ converging a.e. to $u$. Applying Theorem \ref{th:cfrs}, $(u_n)$ is bounded in $C^\alpha(\super{B(x_0,r_0/2}))$ and so $u\in C^\alpha(\super{B(x_0,r_0/2}))$ and then $u\in C^{2,\alpha}(\super{B(x_0,r_0/4}))$ by standard elliptic regularity. By a standard covering argument, we deduce that $u\in C^{2,\alpha}(\super{B_{1/2}})$. 
\hfill\qed

\medskip

As already observed in Remark \ref{Rem1}, Theorems \ref{theorem:liouville0} and \ref{theorem:liouville0b} fail for solutions with positive and finite Morse index. Nevertheless, we can prove radial symmetry and sharp asymptotic behavior at infinity of such solutions. More precisely, we have the following two results.

\medskip

\begin{theorem} \label{theorem:asympt}
Let $u\in C^2(\R^N)$ be a solution of 
	\begin{equation}
	-\Delta u = f(u)\quad\text {in $\R^N$},
	\end{equation}
which is stable outside the ball $B_1$ and bounded below. Assume that $f\geq 0$ is locally lipschitz continuous and $1\le N\le 10$. Then, 

\begin{itemize}

\item [i)] if $ N=1,2$, then $u$ is constant.

\item [ii)] if $ 3 \leq N \leq 9$,  there exists a constant $C>0$ depending on $u$ and $N$ only such that
	\begin{equation}\label{limite}
\vert u(x)-\inf_{\R^N} u\vert\le C\vert x\vert^{-\frac N2 -\sqrt{N-1}+2},\qquad \text{for all $x\in\R^N$,}
	\end{equation} 
	$\vert\nabla u\vert\in L^2(\R^N)$ and 
$$
\int_{\R^N\setminus B_{R}}\vert\nabla u\vert^2\le C R^{-2(\sqrt{N-1}-1)},\qquad \text{for all $R>2$.}
$$

\item [iii)] if $ N = 10$, for any $ \varepsilon >0$ (small enough) there exists a constant $C_{\varepsilon}>0$ depending on $u$, $N$ and $\varepsilon$ only such that 
\begin{equation}\label{limite10}
\vert u(x)-\inf_{\R^N} u\vert\leq C_{\varepsilon} \vert x\vert^{-\frac N2 -\sqrt{N-1}+2 +  \varepsilon},\qquad \text{for all $x\in\R^N$,}
\end{equation} 
$\vert\nabla u\vert\in L^2(\R^N)$ and 
$$
\int_{\R^N\setminus B_{R}}\vert\nabla u\vert^2\le C_{\varepsilon} R^{-2(\sqrt{N-1}-1) + \varepsilon },\qquad \text{for all $R>2$.}
$$
\end{itemize}

Furthermore, $f(\inf_{\R^N} u) =0$ and $ (u - \inf_{\R^N} u) \in L^{2^*}(\R^N)$.
\end{theorem}

\begin{remark}\label{remarksoac}\hfill

\begin{itemize}
	
\item Clearly, Theorem \ref{theorem:asympt} remains true for solutions stable outside a compact set.
\item Theorem \ref{theorem:asympt} complements some results established in \cite{DuFa, Fa3}. It was already known to hold for any $N\ge 3$ in the particular case where $u$ is bounded and radial, see \cite{villegas2}. 
\item Theorem \ref{theorem:asympt} fails if $N\ge12$ for nonradial solutions. Indeed, if $n\ge 11$, there exists a nontrivial bounded stable radial solution $v$ for the nonlinearity $f(v)=v^p$, $p$ large, see \cite{fa2}. Then, the function $u(x',x_N)=v(x')$ for $(x',x_N)\in\R^n\times\R$ is a bounded stable solution in dimension $N=n+1\ge12$. In addition, $u \not \to 0$ as $ \vert x \vert \to \infty $ and $\nabla u\not\in L^2(\R^N)$. 
\item Whether Theorem \ref{theorem:asympt} holds in dimension $N=11$ is an open problem.
\end{itemize}
\end{remark}

\

\begin{theorem} \label{theorem:symmetry}
Let $u\in C^2(\R^N)$ be a solution of \eqref{equazione}
which is stable outside a compact set and bounded below. Assume that $f\in C^{3,1}_{loc}(\R)$, $f\geq 0$ and $N \leq 10$. Then, $u$ is radial symmetric about some point and radially decreasing (and radially strictly decreasing if $u$ is not a constant).
\end{theorem}

\begin{remark}\hfill

\begin{itemize}
\item The assumption $f \in C^{3,1}_{loc}(\R)$ is most likely technical. It is needed only at every zero of $f$. As follows from its proof, the theorem remains true if $f\in C^{1,1}_{loc}(\R)$ and $5\le N\le10$ and if $f\in C^{2,1}_{loc}(\R)$ and $N=4$. 

\item Whether the result remains true in dimension $N=11$ is an open problem, even for the nonlinearity $f(u)=u^p$, $p$ large.  However, the cylindrical solution mentioned in Remark \ref{remarksoac} shows that Theorem \ref{theorem:symmetry} fails in dimension $N\ge12$. Working similarly with the exponential nonlinearity, Theorem \ref{theorem:symmetry} also fails in dimension $N\ge11$ under the weaker lower-bound \eqref{log}.
\item Thannks to the result of \cite{villegas2}, \eqref{limite10} holds true for $\epsilon=0$, that is, optimal asymptotic bounds also hold in dimension $N=10$.
\end{itemize}
\end{remark}

\medskip

\section{Stable solutions}

Let us  prove Theorem \ref{theorem:liouville0}.

\medskip

\noindent {\bf Proof of Theorem \ref{theorem:liouville0}. }

We start by noting  that, if $u$ is a stable solution of \eqref{equazione} in $\R^n, n \geq1$, then the function $v_k : =v_k(x_1,...,x_n,...,x_{n+k}):=u(x_1,...,x_n)$ is a stable solution of \eqref{equazione} in $\R^{n+k}$ for any $k \geq1$. Therefore, it is enough to prove the claim for solutions $u$ in dimension $ N=10$.

Since $u$ is superharmonic and bounded below, its spherical average $\fint_{\partial B_R}u$ decreases to a limit $\l\in\R$, as $R\to+\infty$. Replacing $u$ by $u-\l$, we may assume from here on that
\begin{equation}\label{lim0}
\fint_{\partial B_R}u\searrow 0,\qquad\text{as $R\to+\infty$}
\end{equation}
(note that $ u- \l $ is a bounded below stable solution of \eqref{equazione} with $f$ replaced by the non negative nonlinearity $f(\cdot + \l)$.)

 We begin by proving the following lemma (which holds in any dimension):
\begin{lemma} \label{lemma1}
\begin{equation}\label{smallo}
\int_{B_R}\vert x\vert^{2-N}\vert\nabla u\vert^2 = o(\ln R),\qquad\text{as $R\to+\infty$.}
\end{equation}
\end{lemma}

\proof
According to Proposition 2.5 in \cite{cfrs}, the following $H^1$-bound holds true
$$
\Vert \nabla u\Vert_{L^2(B_1)} \le C \Vert u\Vert_{L^1(B_2)}
$$
Applying the above to $u_R(x)=u(Rx)$ (as we may) yields

\begin{equation}\label{h1bound}
\int_{B_R}\vert\nabla u\vert^2 \le C R^{N-2} \Big (
 \fint_{{B_{2R}}}u \Big )^2    . 
\end{equation}

Here, we have used that $\int_{\partial B_r}u\ge 0$, thanks to \eqref{lim0}.
Using polar coordinates and integration by parts, we find
$$
\int_{B_{R}\setminus B_1}\vert x\vert^{2-N}\vert\nabla u\vert^2 = \int^{R}_1 r^{2-N}\int_{\partial B_r}\vert\nabla u\vert^2 =\left[  r^{2-N}\int_{B_r}\vert\nabla u\vert^2\right]_{r=1}^{r=R}  -  (2-N)
\int_1^{R} r^{1-N}\int_{B_r}\vert\nabla u\vert^2$$
Let us inspect each term in the right-hand side. Since $u$ is superharmonic, $\fint_{B_r}u\le u(0)$. By \eqref{h1bound}, we deduce that the first term is bounded.
For the second term, either $\int_{B_r}u$ is bounded and so 
	$$
	\lim_{r\to+\infty} r^{-N}\int_{B_r}u=0
	$$
	or using L'H\^opital's rule and \eqref{lim0} we have, once again,
$$
\lim_{r\to+\infty} r^{-N}\int_{B_r}u= \lim_{r\to+\infty} N^{-1} r^{1-N}\int_{\partial B_r}u=0
$$
Hence, by \eqref{h1bound} and the above
$$
r^{1-N}\int_{B_r}\vert\nabla u\vert^2 \le C \frac 1r 
\Big ( \fint_{{B_{2r}}} u \Big )^2 = o(1/r), \qquad\text{as $r\to+\infty$.}
$$
Hence,
$$
\int_1^{R} r^{1-N}\int_{B_r}\vert\nabla u\vert^2 = o\left(\int_1^{R} \frac1r\right)=o (\ln R),
$$
and the lemma follows since $\int_{B_1}\vert x\vert^{2-N}\vert\nabla u\vert^2$ is bounded.
\hfill\qed

\

Next, we use inequality\footnote{This inequality is derived by computing the second variation of energy along dilations, similarly yet differently from what can be done in Pohozaev's identity or certain monotonicity formulae.} (2.2) in \cite{cfrs}, which, up to rescaling, reads as follows: for all $\zeta\in C^{0,1}_c(\R^N)$, if $2\le N\le 10$, there holds
\begin{align*}
0\le& \int -2\vert x\vert^{2-N}\vert\nabla u\vert^2\zeta(x\cdot\nabla\zeta)+
\int 4\vert x\vert^{2-N}(x\cdot\nabla u) \zeta \nabla u\cdot\nabla\zeta+\\
&\int(2-N)\vert x\vert^{-N}\vert x\cdot \nabla u\vert^2\zeta(x\cdot\nabla\zeta)+
\int \vert x\vert^{2-N}\vert x\cdot\nabla u\vert^2\vert\nabla \zeta\vert^2
\end{align*}
In particular, if $\zeta$ is radial and $r=\vert x\vert$, the above inequality reduces to
$$
2\int r^{2-N}\vert\nabla u\vert^2\zeta r\zeta' \le \int r^{2-N}\left(\frac{\partial u}{\partial r}\right)^2 r\zeta'\left\{(6-N)\zeta+r\zeta'\right\}
$$
Choose $\zeta$ as follows. Given $R>R_1>2$, $r\in\R_+$,
$$
\zeta(r)=
\left\{
\begin{aligned}
r^4 \qquad\text{in $[0,R_1)$,}\\
R_1^4\qquad\text{in $[R_1,R)$}\\
R_1^4\frac{\ln\frac r{R^2}}{\ln\frac1R}\qquad\text{in $[R,R^2)$,}\\
0\qquad\text{otherwise}
\end{aligned}
\right.
$$
Since $ N=10$, we have that 
$$
r\zeta'\left\{(6-N)\zeta+r\zeta'\right\} = 0\qquad\text{in $[0,R]$,}
$$
which leads us to an inequality of the form
$$
8\int_{B_{R_1}}\vert\nabla u\vert^2 \le CR_1^8 \int_{B_{R^2}\setminus B_R} \vert x\vert^{2-N}\vert\nabla u\vert^2 \left((\ln R)^{-1}+(\ln R)^{-2}\right)
$$
Applying Lemma \ref{lemma1} and letting $R\to+\infty$, we arrive at
$$
\int_{B_{R_1}}\vert\nabla u\vert^2 =0.
$$
Letting $R_1\to+\infty$, we deduce that $u$ is constant.\hfill\qed

\

\noindent {\bf Proof of Theorem \ref{th:liouville1}. } Since $f(0)\geq 0$ and $u$ is bounded, we have that $ \frac{\partial u}{\partial x_N} >0$ on $\R^N_+$ (see \cite{BCN}, \cite{Dancer}). Therefore $ \frac{\partial u}{\partial x_N} $ is a positive solution of the linearized equation $-\Delta w  - f'(u)w = 0$ on $\R^N_+$ and so $u$ is a stable solution of  $-\Delta w = f(w) $ on $\R^N_+$ (see for instance \cite{FisSchoen, MossPi}).

The boundedness of $u$, standard elliptic estimates and the monotonicity of $u$ with respect to the variable $x_N$, imply that the function
\begin{equation*}
v(x_1, . . . , x_{N-1}) := \lim_{x_N\to+\infty} u(x)
\end{equation*}
satisfies
\begin{equation}\label{pbinteroN-1}
\begin{cases}
    v \in C^2(\R^{N-1}) \cap L^{\infty}(\R^{N-1}),\\ 
	\, -\Delta v = f(v)& \text{in $\R^{N-1}$},\\
	\quad \quad v \geq 0 & \text{in $\R^{N-1}$}.
\end{cases}
\end{equation}
In addition $v$ is a stable solution of the above problem (see for  instance \cite{BCN, AC}). Here we have used the continuity of $f^{'}$. 

If $f$ satisfies 1. we can apply Theorem \ref{th:cfrs} to \eqref{pbinteroN-1} to infer that $v \equiv c = const.$ Here we have used that $ N-1 \le 9$. The equation then yields $ f(c) = 0$.  

When 2. is in force, we observe that $u(x) \le z$ for any $x \in \R^N_+ ,$ thanks to Lemma 2.4 in \cite{FVarma}.  The latter and the definition of $v$ imply that $v(x) \le z$ for any $x \in \R^N_+$. Therefore $f(v(x)) \geq 0$ for any $x \in \R^N_+$ and $v$ is a bounded stable solution of 
$-\Delta w = g(w)$ in $\R^{N-1}$ with $ g := f {\bf 1}_{[0,z]} \in C^{0,1}(\R,\R)$, $ g \geq 0$. Another application of Theorem \ref{th:cfrs} then yields $v \equiv c = const.$ and so $ f(c) = 0$. 

In both cases we proved that $v \equiv c = const.$ and $ f(c) = 0$. 
To conclude we observe that $ \frac{\partial u}{\partial x_N} >0$ on $\R^N_+$ implies $\sup_{\R^N_+} u = c$ and so $f(\sup_{\R^N_+} u) =0$ . The one-dimensional symmetry of $u$ is then a consequence of Theorem 1 in \cite{BCNMagenes}. 
\hfill\qed

\

\noindent {\bf Proof of Theorem \ref{th:liouville2}. }
By proceeding as in the first part of the proof of Theorem 2 in \cite{dancer2} we get that the even extension of $u$ to $\R^N$ is a bounded below stable solution of \eqref{equazione} in $\R^N, N \leq 10.$ This function must be constant by Theorem \ref{theorem:liouville0}. The latter implies the desired conclusion.
\hfill\qed

\

\noindent {\bf Proof of Theorem \ref{th:liouville3}. }
Assume by contradiction that $u>0$ in $\Omega$. 
According to \cite{BCN3} (see also \cite{el} for a prior result concerning smooth epigraphs)
$u$ is monotone. Therefore $u$ is a stable solution and working as in the proof of Theorem \ref{th:liouville1}, we deduce that
$$v(x_1, . . . , x_{N-1}) := \lim_{x_N\to+\infty} u(x)$$
is a positive stable solution of $ - \Delta u = f(u) $ in $ \R^{N-1}$, with $ N-1 \leq 10$, and so it must be constant by Theorem \ref{theorem:liouville0}. But this is in contradiction with $f(t)>0$ for $t>0$. Hence $u$ must vanish somewhere in $ \Omega$ and so it must be identically zero by the strong maximum principle. The latter then implies that $f(0)=0$.
\hfill\qed

\section{Solutions which are stable outside a compact set}

\noindent {\bf Proof of Theorem \ref{theorem:asympt}. } 
If $ N \leq 2$ the solution u must be constant since it is bounded and superharmonic. We may therefore suppose that $N \geq 3$. 
	
We begin by adapting Lemma \ref{lemma1} to solutions which are stable outside a compact set.
\begin{lemma} \label{lemma2} Assume that $N\ge3$ and that $u$ is a solution which is stable outside a compact set. Then, 
\begin{equation}\label{smallo2}
\int_{B_R}\vert x\vert^{2-N}\vert\nabla u\vert^2 = o(\ln R),\qquad\text{as $R\to+\infty$.}
\end{equation}
\end{lemma}
\proof Assume without loss of generality	 that \eqref{lim0} holds and that $R=2^n$ for some $n\in\N^*$. Then,
$$
\int_{B_{R}\setminus B_2} \vert\nabla u\vert^2=\int_{B_{2^n}\setminus B_2} \vert\nabla u\vert^2 = \sum_{k=2}^n \int_{B_{2^k}\setminus B_{2^{k-1}}}\vert\nabla u\vert^2.
$$
Given $k\in\{2,\dots,n\}$, let $v(x)=u(2^{k-1}x)$, $x\in\R^N\setminus B_{2^{-(k-1)}}$. Then, $v$ is stable outside the ball of radius $B_{2^{-(k-1)}}$ and
$$
\int_{B_{2^k}\setminus B_{2^{k-1}}}\vert\nabla u\vert^2 =2^{(k-1)(N-2)}\int_{B_2\setminus B_1}\vert\nabla v\vert^2
$$
The annulus $B_2\setminus B_1$ can be covered by finitely many balls of radius $\frac12$. On each of these balls, $v$ is stable. Applying Proposition 2.5 in \cite{cfrs}, we deduce that
$$
\int_{B_2\setminus B_1}\vert\nabla v\vert^2\le C\int_{B_{4}\setminus B_{1/2}}\vert v\vert = C\fint_{B_{2^{k+1}}\setminus B_{2^{k-2}}} u =o(1)\quad\text{as $k\to+\infty$},
$$
where we used \eqref{lim0} for the last equality. So,
\begin{equation}\label{annu}
\int_{B_{R}\setminus B_2} \vert\nabla u\vert^2 \le C\sum_{k=2}^n o(2^{(k-1)(N-2)})= o(R^{N-2})\quad\text{as $R=2^n\to+\infty$.}
\end{equation}
Integrate by parts $\int_{B_{R}\setminus B_1}\vert x\vert^{2-N}\vert\nabla u\vert^2$ exactly as in the proof of Lemma \ref{lemma1}
and \eqref{smallo2} follows.\hfill\qed

\medskip

Next, we prove that $\nabla u\in L^2(\R^N)$. Assume that $3\le N\le 9$.  According to Lemma 2.1, inequality (2.2) in \cite{cfrs}, for all $\zeta\in C^{0,1}_c(\R^N)$, $\zeta$ radial with support outside the ball of radius 1, there holds
\begin{multline}\label{pohozaev}
\frac{(N-2)(10-N)}{4} \int_{\R^N} r^{2-N}\left(\frac{\partial u}{\partial r}\right)^2\zeta^2 \le   \\
\le -2 \int_{\R^N} r^{2-N}\vert\nabla u\vert^2\zeta r\zeta' + \int_{\R^N} r^{2-N}\left(\frac{\partial u}{\partial r}\right)^2 r\zeta'\left\{(6-N)\zeta+r\zeta'\right\}
\end{multline}
By analogy with a strategy found in \cite{villegas2}, we choose $\zeta$ as follows. Fix $R_2>R>2$ and 
\begin{equation}\label{alfa}
\alpha = \frac N2+ \sqrt{N-1} -2.
\end{equation}
Given  $r\in\R_+$, let
$$
\zeta(r)=
\left\{
\begin{aligned}
2^\alpha(r-1) \qquad\text{in $[1,2)$,}\\
r^\alpha\qquad\text{in $[2,R)$}\\
R^{\alpha}\qquad\text{in $[R,R_2)$}\\
R^{\alpha}\frac{\ln\frac r{R_2^2}}{\ln\frac1{R_2}}\qquad\text{in $[R_2,R_2^2)$,}\\
0\qquad\text{otherwise}
\end{aligned}
\right.
$$
All integral terms in \eqref{pohozaev} in the region $[1\le\vert x\vert \le 2]$ are controlled by a constant $C$ depending on $u$ and $N$ only. For the region $[2\le\vert x\vert\le R]$, our choice of $\alpha$ leads to the cancellation of all terms involving the radial derivative of $u$, so that all remains is the {\it negative} term (in the right-hand side)
$$
-2\alpha\int_{B_R\setminus B_2} r^{2-N+2\alpha}\vert\nabla_T u\vert^2,
$$
where we denoted the tangential part of the gradient by $\nabla_T u = \nabla u -\frac{\partial u}{\partial r}e_r$, $e_r=x/\vert x\vert$.
In the region $[R,R_2)$, all terms disappear except the left-hand side:
$$
R^{2\alpha}\frac{(N-2)(10-N)}{4}\int_{B_{R_2}\setminus B_R} r^{2-N}\left(\frac{\partial u}{\partial r}\right)^2
$$
Finally, in the region $[R_2\le\vert x\vert\le R_2^2]$, all terms can be bounded above by a constant $C=C(R,N)$ times 
$$
\left(\int_{B_{R_2^2}\setminus B_{R_2}} \vert x\vert^{2-N}\vert\nabla u\vert^2 \right)\left((\ln R_2)^{-1}+(\ln R_2)^{-2}\right). 
$$
This quantity converges to $0$ as $R_2\to+\infty$, thanks to Lemma \ref{lemma2}. So, in the limit $R_2\to+\infty$, inequality \eqref{pohozaev} reduces to
$$
\frac{(N-2)(10-N)}{4}R^{2\alpha}\int_{\R^N\setminus B_R} \vert x\vert ^{2-N}\left(\frac{\partial u}{\partial r}\right)^2  + 2\alpha\int_{B_R\setminus B_2} \vert x\vert^{2-N+2\alpha}\vert\nabla_T u\vert^2\le C.
$$
The above inequality being true for all $R>2$, we readily deduce that for $3\le N\le 9$,
\begin{equation}\label{annulus}
\int_{B_{2R}\setminus B_R}\vert\nabla u\vert^2\le C R^{-(2-N+2\alpha)}
\end{equation}
When $N=10$, we only have 
$$
\int_{B_{2R}\setminus B_R}\vert\nabla_T u\vert^2\le C R^{-(2-N+2\alpha)}
$$
However, replacing $\alpha$ by $\alpha_\epsilon=\alpha-\epsilon$, $\epsilon>0$ small, in the definition of $\zeta$ leads to an inequality of the form
$$
\int_{B_{2R}\setminus B_R}\vert\nabla u\vert^2\le C_\epsilon R^{-(2-N+2\alpha_\epsilon)},
$$
Note that $2-N+2\alpha=2(\sqrt{N-1}-1)>0$ for $N\ge3$. So, applying this inequality with $2^k R$, $k\in \N$, in place of $R$ and summing over $k$ implies that $\nabla u\in L^2(\R^N)$ and 
\begin{equation}\label{grad}
\int_{\R^N\setminus B_R}\vert\nabla u\vert^2\le C R^{-2(\sqrt{N-1}-1)},
\end{equation}
if $3\le N\le 9$ (and $\nabla u\in L^2(\R^N)$ and $\int_{\R^N\setminus B_R}\vert\nabla u\vert^2\le C_{\varepsilon} R^{-2(\sqrt{N-1}-1) + \epsilon} \, $ if $ N=10$).

Next, fix a point $x\in \R^N\setminus \super{B_2}$ and let $R=\vert x\vert/2$. Take another point $y\in \partial B(x,r)$, $r\le R$, and apply the fundamental theorem of calculus:
$$
\vert u(x)-u(y)\vert\le \int_0^1\left\vert\frac d{dt}u(x+t(y-x))\right\vert dt\le r\int_0^1\vert\nabla u(x+t(y-x))\vert\;dt
$$
Integrating over $\partial B(x,r)$, we deduce that
$$
\int_{\partial B(x,r)}\vert u(x)-u(y)\vert\;d\sigma(y)\le \int_{B(x,r)}\vert\nabla u\vert\le \int_{B(x,R)}\vert\nabla u\vert
$$
Integrating once more in $r\in(0,R)$, using Cauchy-Schwarz and applying at last \eqref{grad}, for $ 3 \leq N \leq 9$ we get 
\begin{multline}\label{campa}
\fint_{B(x,R)}\vert u(x)-u(y)\vert\;dy \le CR^{1-N}\int_{B(x,R)}\vert\nabla u\vert\le CR^{1-N/2}\left(\int_{B(x,R)}\vert\nabla u\vert^2\right)^{\frac12}\\
\le CR^{1-N/2}\left(\int_{B(0,3R)\setminus B(0,R)}\vert\nabla u\vert^2\right)^{\frac12}
\le CR^{-\frac N2 -\sqrt{N-1}+2},
\end{multline}
(for $ N=10$ the latter is replaced by $ C_{\varepsilon} R^{-2(\sqrt{N-1}-1) + \epsilon}$).

We may draw two conclusions from the above inequality. Firstly, 
$$\lim_{\vert x\vert\to+\infty}u(x)=0.$$ 
Recalling that $u$ is superharmonic, we have indeed  
\begin{align*}
u(x)\le \fint_{B(x,R)} \vert u(x)-u(y)\vert\;dy +\fint_{B(x,R)}u&\le CR^{-\frac N2 -\sqrt{N-1}+2}
+ C\fint_{B(0,3R)}u(y)\;dy\\
&\le CR^{-\frac N2 -\sqrt{N-1}+2}+Cu(0)
\le C.
\end{align*}

So, $u$ is bounded. By elliptic regularity, so is $\nabla u$. 
	We also know that $\nabla u\in L^2(\R^N)$, whence $u-\xi \in L^{2^*}(\R^N)$, for some $ \xi \in \R$, thanks to Theorem 1.78 of \cite{MZ}. The proof of this result also implies that $ \xi = \lim_{j\to\infty} \fint_{B(0, 2^j)} \, u $ and so $ \xi =0$ by \eqref{lim0} (as already seen in the proof of Lemma \ref{lemma1}.) Therefore, $u\in L^{2^*}(\R^N)$ and thus $u^{2*}$ is integrable and globally Lipschitz on $ \R^N$. This clearly implies that $\lim_{\vert x\vert\to+\infty}u(x)=0.$

Secondly, similarly to \eqref{campa}, we have for $3\le N\le 9$
\begin{align*}
\fint_{B(2x,R)}\vert u(2x)-u(y)\vert\;dy &\le CR^{-\frac N2 -\sqrt{N-1}+2}\quad\text{and}\\
\fint_{B\left(\frac32 x,2R\right)}\left\vert u\left(\frac32 x\right)-u(y)\right\vert\;dy &\le CR^{-\frac N2 -\sqrt{N-1}+2}.
\end{align*}
Using the notation $u_{z,r}=\fint_{B(z,r)} u$ for the average of $u$ on a given ball $B(z,r)$, it follows that
\begin{align*}
\vert u(x) - u(2x)\vert &\le \vert u(x) - u_{x,R}\vert+ \left\vert u_{x,R} - u\left(\frac32 x\right)\right\vert+ \left\vert u\left(\frac32 x\right)- u_{2x,R}\right\vert + \vert  u_{2x,R}-u(2x) \vert\\
&\le  CR^{-\frac N2 -\sqrt{N-1}+2} + \fint_{B(x,R)}  \left\vert u\left(\frac32 x\right) - u(y)\right\vert dy +  \fint_{B(2x,R)}  \left\vert u\left(\frac32 x\right) - u(y)\right\vert dy\\
&\le CR^{-\frac N2 -\sqrt{N-1}+2} + C \fint_{B\left(\frac32 x,2R\right)}  \left\vert u\left(\frac32 x\right) - u(y)\right\vert dy \\
&\le C\vert x\vert^{-\frac N2 -\sqrt{N-1}+2}
\end{align*}
Applying the above inequality to $2^k x$, $k\in\N$, in place of $x$ and summing over $k$, we deduce that the sequence $(u_k)$ defined by $u_k(x)=u(2^kx)$, 
	converges in $C^0_{loc}(\R^N \setminus \overline{B_2})$ to a limit $v$ as $k\to+\infty$ and
	$$
	\vert u(x)-v(x)\vert\le C\vert x\vert^{-\frac N2 -\sqrt{N-1}+2} \qquad {\text {if}} \quad 3 \leq N \leq 9,
	$$
	(resp. $ \vert u(x)-v(x)\vert\le C_{\varepsilon} \vert x\vert^{-\frac N2 -\sqrt{N-1}+2 + \varepsilon} \, $ if $ N=10$). 
	Since $\lim_{\vert x\vert\to+\infty}u(x)=0$ we necessarily have that $ v \equiv 0$ and so 
	$$
	\vert u(x)\vert\le C\vert x\vert^{-\frac N2 -\sqrt{N-1}+2} \qquad {\text {if}} \quad 3 \leq N \leq 9,
	$$
	(resp. $ \vert u(x)\vert\le C_{\varepsilon} \vert x\vert^{-\frac N2 -\sqrt{N-1}+2 + \varepsilon} \, $ if $ N=10$). 

In order to establish \eqref{limite} and \eqref{limite10}, it remains to prove that $\inf_{\R^N}u=0$. We consider a sequence  $(x_n)$ such that $u(x_n) \to \inf_{\R^N}u$. Then, either $\vert x_n \vert \to +\infty$ or $(x_n)$ posseses a bounded subsequence (still called $(x_n))$ such that $ x_n \to { \overline x}$.  In the first case  we clearly have $\inf_{\R^N}u=0$,  while in the second case we get that $u({ \overline x}) =\inf_{\R^N}u$ and so the nonnegative subharmonic function $v := u - u({ \overline x})$ vanishes at ${ \overline x}$. Therefore $v$ must be zero by the strong maximum principle which, in turn, yields that $u \equiv \inf_{\R^N}u=\lim_{\vert x\vert\to+\infty}u(x)=0$.
\hfill\qed

\noindent{\bf Proof of Theorem \ref{theorem:symmetry}.}

If $ N \leq 2$, $u$ is constant by Theorem  \ref{theorem:asympt} and so we are done. 
When $ N\geq 3$ we can assume that $u$ is not constant  (otherwise the claim is trivially true). As before, up to replacing $u$ by $u - \inf_{\R^N} u$, we may and do suppose that $ \inf_{\R^N} u =0$ and so $ u > 0$ in $\R^N$ by the strong maximum principle. Also, thanks to the asymptotics of Theorem \ref{theorem:asympt} we have  

\begin{equation}\label{int-u}
\begin{cases}
u \in L^{\frac{N}{2}}(\R^N) \quad \,\,   \text{if} \qquad N \geq 5,\\
u^2 \in  L^{\frac{N}{2}}(\R^N) \quad \text{if} \qquad N =4 ,\\
u^3 \in  L^{\frac{N}{2}}(\R^N) \quad \text{if} \qquad N=3.
\end{cases}
\end{equation}

Moreover, since $ f \geq 0$ by assumption and $f(0)=0$ by Theorem \ref{theorem:asympt}, we must have $f'(0)=0$ and $f''(0) \geq0$. We also observe that in dimension $N=3,4$ we must have  $ f''(0)=f'''(0) =0$.  Indeed,  if $f''(0)> 0$, then $ \liminf_{t\to 0^+} \frac{f(t)} {t^{\frac{N}{N-2}}} \in (0, +\infty]$ and Theorem 3.5 of \cite{DambMiti} implies that $u$ must be identically zero, in contradiction with $u>0$. Thus,  $f''(0) = 0$ and so $f'''(0) =0$ (again by $f \geq 0$). In particular we have 

\begin{equation}\label{comp-f}
\forall \, t \in [0, \max_{\R^N}u] \qquad \qquad \vert f'(t) \vert \leq 
\begin{cases}
C t \quad \,\, \text{if} \qquad N \geq 5,\\
C t^2 \quad \text{if} \qquad N =4 ,\\
C t^3 \quad \text{if} \qquad N=3,
\end{cases}
\end{equation}
where $ C>0$ is a constant depending only on $f$. 

We are now ready to apply the moving planes method. Given $\lambda\in\R$, set $\Sigma_\lambda=\{(x_1,x')\in\R^N\;:\; x_1<\lambda\}$ and for $x\in\Sigma_\lambda$, let as usual $x^\lambda=(2\lambda-x_1,x')$, $u_\lambda(x)=u(x^\lambda)$ and $w_\lambda=u_\lambda-u$. Since $u>0$ and $ \lim_{\vert x \vert \to +\infty} u(x) =0$ we may and do suppose that (up to translation) $u(0) = \max_{\R^N} u>0$ and then we prove that $u\equiv u_{0}$ in $\Sigma_{0}$.
Since the coordinate axis $x_1$ can be chosen arbitrarily, we will then conclude that $u$ is radially symmetric about the origin.  That $u$ is radially strictly decreasing is also a standard consequence of the moving planes procedure and the strong maximum principle. This will provide the desired result. 

We claim that $w_\lambda\ge 0$ in $\Sigma_\lambda$ if $\lambda\le -K$ for some $K>0$. Indeed we have that
\begin{equation}\label{eq:w}
-\Delta w_\lambda  = f(u_\lambda)-f(u) = a_\lambda w_\lambda  \qquad\text {in} \quad \Sigma_\lambda,
\end{equation}
where 
\begin{equation}
a_\lambda := \begin{cases}
\frac{f(u_\lambda)-f(u)}{u_\lambda-u} \quad \text{if} \quad u - u_\lambda>0, \\
\qquad 0 \qquad \quad \text{if} \quad u - u_\lambda \leq 0,
\end{cases}
\end{equation}
belongs to $L^{N/2}(\Sigma_\lambda)$. To see this we observe that, by the mean value theorem, for every $ x \in \Sigma_\lambda$ we have  $\vert a_\lambda (x) \vert = \vert f'(\xi(x,\lambda)) \vert {\textbf 1 }_{\{u - u_\lambda >0\}}(x)$, where $\xi(x,\lambda) \in (u_\lambda (x) , u(x))$. Combining the latter with \eqref{comp-f} we get
\begin{equation}\label{int-a}
\vert a_\lambda \vert \leq 
\begin{cases}
C u\quad \,\, \text{if} \qquad N \geq 5,\\
C u^2 \quad \text{if} \qquad N =4 , \qquad \qquad {\text {on}} \quad \Sigma_{\lambda}\\
C u^3 \quad \text{if} \qquad N=3,
\end{cases}
\end{equation}
where $ C>0$ is a constant depending only on $f$. The latter and \eqref{int-a} imply  $ a_\lambda \in L^{N/2}(\Sigma_\lambda)$ for any $ \lambda$. 
Therefore, given any $ \varepsilon >0$ we can find $K>0$ such that 
 $\Vert a_\lambda \Vert_{L^{N/2} (\Sigma_\lambda)}  <\epsilon$, for any $ \lambda \leq -K.$

Let $\chi_R(x)=\chi(x/R)$ be a standard cut-off function, with $R>1$.  
	Multiplying \eqref{eq:w} by $w_\lambda^-\chi_R^2$, integrating by parts and making use of Sobolev's inequality we find
	$$
	\left(\int_{\Sigma_\lambda} \vert w_\lambda^-\chi_R\vert^{2^*}\right)^{2/2^*}\le C\int_{\Sigma_\lambda}   \vert\nabla[w_\lambda^-\chi_R]\vert^2\le C\left\{ \int_{\Sigma_\lambda} (w_\lambda^-)^2\vert\nabla\chi_R\vert^2+\Vert a_\lambda\Vert_{L^{N/2}(\Sigma_\lambda)} \left(\int_{\Sigma_\lambda} \vert w_\lambda^-\chi_R\vert^{2^*}\right)^{2/2^*}\right\}
	$$
	Letting $R\to+\infty$ and then choosing $\Vert a_\lambda\Vert_{L^{N/2}(\Sigma_\lambda \cap [u>u_\lambda])}<\epsilon$ small enough, we deduce that $(w_\lambda)^-\equiv0$ as claimed.  Here we also used that $(w_\lambda^-)^2 \in L^{\frac{N}{N-2}}(\R^N)$ to show that  $\int_{\Sigma_\lambda} (w_\lambda^-)^2\vert\nabla\chi_R\vert^2 \to 0$ as $R\to\infty$ (recall that $u \in L^{2^*}(\R^N)$ by Theorem \ref{theorem:asympt}).

	Let us now finish the moving planes method by setting
	$$
	\lambda_0:=\sup\{\lambda <0 \, : \,  w_t\ge0 \quad\text{in $\Sigma_t$} \quad \forall t \leq \lambda \, \}.
	$$
	To conclude the proof it is enough to prove that $\lambda_0=0$.  
	We argue by contradiction and suppose that	$\lambda_0 <0.$ 
    By continuity, $w_{\lambda_0}\ge0$ in $\Sigma_{\lambda_0}$. By the strong maximum principle we deduce that either $w_{\lambda_0}>0$ in $\Sigma_{\lambda_0}$ or $w_{\lambda_0} \equiv 0$ in $\Sigma_{\lambda_0}$. The latter is not possible if $ \lambda_0 < 0$ since in this case we would have $ u(2\lambda_0,0') = u(0) = \max_{\R^N} u>0 $ and so $ w_\lambda \equiv 0$ in $ \Sigma_{\lambda} $ for any $ \lambda \in [2\lambda_0, \lambda_0]$ by the strong maximum principle and the Hopf lemma. By repeating (infinitely many times) the previous argument we would find $ u(s,0') = u(0) = \max_{\R^N} u >0 $ for any $ s < \lambda_0$, contradicting thus the assumption $ \lim_{\vert x \vert \to +\infty} u(x) =0$. Therefore $w_{\lambda_0}>0$ in $\Sigma_{\lambda_0}$.
	
Now we achieve a contradiction by proving the existence of $\tau_0>0$ such that for any $0<\tau<\tau_0$ we have $w_{\lambda_0+\tau} \geq 0$ in
$\Sigma_{\lambda_0+\tau}$. 
 
To this end we are going to show that, for every $ \delta>0$ there are $ \tau_0>0$ 
and a compact set $	K \subset \Sigma_{\lambda_0}$  (both depending on $\delta$ and $ \lambda_0$) such that 
\begin{equation}\label{fine-mp}
	\Vert a_\lambda\Vert_{L^{N/2}(\Sigma_\lambda \setminus K)}< \delta, \qquad w_{\lambda}>0 \quad {\text {in}} \,\,  K, \qquad \forall \,  \lambda \in [\lambda_0, \lambda_0 + \tau_0].
\end{equation}
We prove this for $ N \geq 5$ (the case $N=3,4$ is obtained in the same way by using \eqref{int-u} and  \eqref{int-a} for $N=3,4$). First pick a compact set $K \subset \Sigma_{\lambda_0}$ such that $\int_{\Sigma_{\lambda_0 \setminus K}} (Cu)^{\frac{N}{2}} < (\frac{\delta}{2})^{\frac{N}{2}}$ and then fix a $\tau_1= \tau_1(\delta,\lambda_0)>0$ such that $\int_{\Sigma_{\lambda \setminus \lambda_0}} (Cu)^{\frac{N}{2}} < (\frac{\delta}{2})^{\frac{N}{2}}$ for every $ \lambda \in [\lambda_0, \lambda_0 + \tau_1],$ where $C$ is the constant appearing in \eqref{int-a} (this choice is clearly possible in view of \eqref{int-u}). Combining these information with \eqref{int-a} we obtain that $	\Vert a_\lambda\Vert_{L^{N/2}(\Sigma_\lambda \setminus K)}< \delta$ for any $ \lambda \in [\lambda_0, \lambda_0 + \tau_1].$ Then, the uniform continuity of the function $  u(2\lambda-x_1,x^{'}) - u(x)$ on the compact set $K \times [\lambda_0, \lambda_0 + \tau_1]$ and the fact that $w_{\lambda_0}>0$ in $\Sigma_{\lambda_0}$ ensure that $ u_{\lambda_0+\tau} - u >0 $ in $K$ for any $0 \le \tau< \tau_2$, for some $\tau_2= \tau_2(\delta,\lambda_0) < \tau_1$. Hence, \eqref{fine-mp} holds true with $\tau_0 \in (0, \tau_2).$

As before, we multiply the equation \eqref{eq:w} by $w_\lambda^-\chi_R^2$, we  integrate by parts and use Sobolev's inequality to get 
$$
\left(\int_{\Sigma_\lambda} \vert w_\lambda^-\chi_R\vert^{2^*}\right)^{2/2^*}\le C_S^2\int_{\Sigma_\lambda}   \vert\nabla[w_\lambda^-\chi_R]\vert^2\le C_S^2\left\{ \int_{\Sigma_\lambda} (w_\lambda^-)^2\vert\nabla\chi_R\vert^2+\Vert a_\lambda\Vert_{L^{N/2}(\Sigma_\lambda)} \left(\int_{\Sigma_\lambda} \vert w_\lambda^-\chi_R\vert^{2^*}\right)^{2/2^*}\right\}
$$
which, in view of \eqref{fine-mp} with $ \delta = \frac{1}{2 C_S^2}$, gives for every $\lambda \in [\lambda_0, \lambda_0 + \tau_0]$,
$$
\left(\int_{\Sigma_\lambda \setminus K} \vert w_\lambda^-\chi_R\vert^{2^*}\right)^{2/2^*}\le C_S^2\left\{ \int_{\Sigma_\lambda \setminus K} (w_\lambda^-)^2\vert\nabla\chi_R\vert^2+\Vert a_\lambda\Vert_{L^{N/2}(\Sigma_\lambda)} \left(\int_{\Sigma_\lambda \setminus K} \vert w_\lambda^-\chi_R\vert^{2^*}\right)^{2/2^*}\right\}
$$
$$
\le C_S^2 \int_{\Sigma_\lambda \setminus K} (w_\lambda^-)^2\vert\nabla\chi_R\vert^2+
\frac {1}{2}  \left(\int_{\Sigma_\lambda \setminus K} \vert w_\lambda^-\chi_R\vert^{2^*}\right)^{2/2^*}.
$$
Then
$$
\left(\int_{\Sigma_\lambda \setminus K} \vert w_\lambda^-\chi_R\vert^{2^*}\right)^{2/2^*}\le 2 C_S^2 \int_{\Sigma_\lambda \setminus K} (w_\lambda^-)^2\vert\nabla\chi_R\vert^2 \to 0 \qquad {\text {when}} \quad R \to \infty,
$$
since $u \in L^{2^*}(\R^N).$ Therefore, $\int_{\Sigma_\lambda \setminus K} \vert w_\lambda^-\vert^{2^*} =0$ and so $ w_\lambda \geq 0$ in $ \Sigma_\lambda$, for every $\lambda \in [\lambda_0, \lambda_0 + \tau_0]$. This contradicts the definition of $ \lambda_0$ and so $\lambda_0 =0$ which, in turn, yields $ u \leq u_0$ in $ \Sigma_0$.  Now we can apply the same procedure to the function 
$ v(x):= u(-x_1,x') (= u_0(x))$ to find that $ v \le v_0$  in $ \Sigma_0$, i.e., $ u_0 \leq u$ in $ \Sigma_0$. This proves that $ u_0 \equiv u$ on $ \Sigma_0$ and concludes the proof.

\eproof

\bibliographystyle{amsalpha}

\end{document}